\documentclass[11pt,reqno]{amsart}

\usepackage{amssymb,amsthm}
\usepackage{ccaption}
\usepackage{mathrsfs}
\usepackage[margin=30mm]{geometry}
\usepackage{xcolor}
\usepackage{array}

\usepackage{comment}

\newtheorem{thm}{Theorem}[section]

\newtheorem*{mthm}{Main Theorem}

\theoremstyle{definition}
\newtheorem{dfn}[thm]{Definition}
\newtheorem{ex}[thm]{Example}

\makeatletter
	
	\@addtoreset{equation}{section}
 \makeatother

\newcommand{\bs}[1]{{\boldsymbol #1}}

\title[Moduli space of contact instantons]{Moduli spaces of contact instantons
on Sasakian 5-manifolds with transverse Calabi-Yau structures
and orbifold $K3$ surfaces}
\author[T.~Arai]{Tomohiro Arai}
\author[K.~Baba]{Kurando Baba}


\address[T.~Arai]{Department of Mathematics, Graduate School of Science and Technology,
Tokyo University of Science,
2641 Yamazaki, Noda, Chiba, 278-8510, Japan.}
\email{6124503@ed.tus.ac.jp}

\address[K.~Baba]{Faculty of Science and Technology,
Tokyo University of Science,
2641 Yamazaki, Noda, Chiba, 278-8510, Japan.
Research Institute for Science and Technology at Tokyo University of Science,
Division of Joint Research of Geometry and Natural Science}
\email{kurando.baba@rs.tus.ac.jp}

\date{}

\begin{document}

\begin{abstract}
We study anti-self-dual contact instantons on 5-dimensional Sasakian manifolds
with transverse Calabi-Yau structures. 
In this case,
the leaf space is a Calabi-Yau orbifold, and 
the moduli space of irreducible anti-self-dual contact instantons
is a hyper\"ahler manifold.
Using the singularity data of the leaf spaces,
we prove that the transverse Levi-Civita connection gives an irreducible anti-self-dual contact instanton
in the case when the leaf space is one of the 95 orbifold $K3$ surfaces classified by Reid. 
Moreover, we compute explicitly the complex dimension of the corresponding moduli spaces
in all 95 cases.
\end{abstract}

\maketitle
%\tableofcontents

\section{Introduction}

In recent years,
the study of contact instantons
has attracted considerable attention
at the intersection of
higher-dimensional gauge theory
and contact geometry.
Contact instantons provide
a natural extension of
$4$-dimensional instanton theory
to $5$-dimensional contact manifolds,
and are closely related to
supersymmetric gauge theory
in mathematical physics.
K\"all\'en-Zabzine (\cite{KZ})
introduced the notion of contact instantons
in the context of $\mathcal{N}=1$
supersymmetric Yang-Mills theory
and showed that contact instantons arise
from the BPS condition
that minimizes the supersymmetric Yang-Mills action
on $5$-dimensional contact manifolds.
Subsequently,
Baraglia-Hekmati
(\cite{BH})
reformulated contact instantons
in the framework of contact geometry
and rigorously constructed
the transverse analogue of
the BPS equations
on the contact distributions of contact manifolds.
They also constructed the moduli space of
anti-self-dual (ASD) contact instantons
and established fundamental geometric properties
of this moduli space.
In particular,
they showed that its dimension depends only
on the basic cohomology
of the underlying contact manifold.

The aim of the present paper is
to study ASD contact instantons
and their moduli space
on $5$-dimensional Sasakian
manifolds
with transverse Calabi-Yau structures.
Such contact manifolds
provides us a concrete geometric setting
in which Baraglia-Hekmati's framework can
be applied.
To state our main theorem,
we recall the definition of
contact instantons following to \cite{BH}.
Let $M$ be a $5$-dimensional Sasakian manifold
with contact form $\eta$.
Let $P$ be a principal fiber bundle over
$M$ whose structure group is a compact Lie group
$G$.
A connection
$\nabla$ on
$P$
(or on its adjoint bundle
$\mathfrak{g}_{P}$)
is called a \textit{self-dual contact instanton}
(SD contact instanton)
if its curvature $2$-form
$F^{\nabla} \in \Omega^{2}(\mathfrak{g}_{P})$
satisfies
\begin{equation}\label{eqn:SDCI_eq}
  \ast F^{\nabla} = \eta \wedge F^{\nabla},
\end{equation}
and an \textit{anti-self-dual contact instanton}
(ASD contact instanton)
if it satisfies
\begin{equation}\label{eqn:ASDCI_eq}
  \ast F^{\nabla} = -\eta \wedge F^{\nabla},
\end{equation}
where $\ast$ denotes the Hodge star operator on
$M$.
Equations
\eqref{eqn:SDCI_eq}
and
\eqref{eqn:ASDCI_eq}
are called the self-dual
and the anti-self-dual contact instanton equations,
respectively.
It should be noted that
any ASD contact instanton
is automatically a Yang-Mills connection,
while the same property does not necessarily hold
for SD contact instantons.

Among Sasakian manifolds,
compact Sasakian manifolds
with transverse Calabi-Yau structures
provide a natural setting
for the study of ASD contact instantons.
On such manifolds,
one can construct the moduli spaces
of gauge equivalence classes
of ASD contact instantons.
Baraglia-Hekmati
\cite{BH}
introduced the basic deformation complex
associated with the space
$\Omega^{*}_{B}(M)\otimes\mathfrak{g}_{P}$
of $\mathfrak{g}_{P}$-valued basic differential forms
and showed that the moduli space
$\mathscr{M}^{*}$
of irreducible ASD contact instantons
is a smooth manifold
(\cite[Corollary 4.4]{BH}).
They also constructed the structure of a hyperk\"ahler manifold
on
$\mathscr{M}^{*}$,
which is induced from the transverse Calabi-Yau structure
of $M$ (\cite[Proposition 4.9]{BH}).

Let $M$ be a Sasakian manifold
with a transverse Calabi-Yau structure.
Let $\mathcal{F}_{\xi}$ denote
the characteristic foliation on
$M$,
which is generated by the Reeb vector field
$\xi$ on $M$.
By decomposing the Levi-Civita connection
$\nabla^{g}$
of
$M$
into its horizontal and vertical components
along $\mathcal{F}_{\xi}$,
one obtains the transverse Levi-Civita connection $\nabla^{*}$
on the contact distribution $H = \operatorname{Ker}(\eta)$.
It is shown in \cite[Proposition~2.11]{BH}
that $\nabla^{*}$ is an ASD contact instanton
on the adjoint bundle
$\mathfrak{g}_{P}=\Lambda^{-}H^{*}$
of a certain
principal $\mathrm{SO}(3)$-bundle over
$M$
(see \eqref{eqn:bundle-decomp} for the definition of $\Lambda^{-}H^{*}$).
By restricting attention to
the class of Sasakian manifolds
whose leaf spaces $X=M/\mathcal{F}_{\xi}$
are 
codimension one,
orbifold $K3$ surfaces
in the weighted projective space
$\mathbb{CP}(w_{0},w_{1},w_{2},w_{4})$
which was classified by Reid \cite{Reid1, Reid2}
(often referred to as
\textit{Reid orbifold $K3$ surfaces}
in this paper).
Then we can investigate
the structure of the moduli space
in greater detail by using the singularity data of $X$.
In fact,
we explicitly determine
the irreducibility
of the transverse Levi-Civita connection
$\nabla^{*}$
and the complex dimension
of the moduli space
$\mathscr{M}^{*}$.

\begin{mthm}\label{thm:main_result}
Let $X$ be a Reid orbifold $K3$ surface
and $M$ be a compact simply-connected
$5$-dimensional Sasakian manifold
with a transverse Calabi-Yau structure
whose leaf space is $X$.
Then the transverse Levi-Civita connection
$\nabla^{*}$
is an irreducible ASD contact instanton,
and the moduli space
$\mathscr{M}^{*}$
of irreducible ASD contact instantons
on $\mathfrak{g}_{P} = \Lambda^{-}H^{*}$
is a non-empty hyperk\"ahler manifold.
The complex dimensions of $\mathscr{M}^{*}$
for all $95$ Reid orbifold $K3$ surfaces
are listed in Table \ref{table:main_result}.
\end{mthm}

We remark that
for any Reid orbifold $K3$ surface
$X$,
there exists a compact $5$-dimensional Sasakian manifold
with a transverse Calabi-Yau structure
whose leaf space is $X$
(see \cite{BG}).
The main theorem extends \cite[Example 5.6]{BH}
to all $95$ families of Reid orbifold $K3$ surfaces
$X$ as the leaf space of
Sasakian manifold
with transverse Calabi-Yau structures.
It should be noted that
$X$ labeled $79$
in Table \ref{table:main_result}
corresponds to the $0$-dimensional
moduli space $\mathscr{M}^{*}$.
This means that the corresponding
transverse Levi-Civita connection
is rigid.

The organization of this paper is as follows.
In Section
\ref{sec:Pre},
we review basic facts
on Sasakian manifolds
with transverse Calabi-Yau structures
and contact instantons,
following to \cite{BG, BH}.
Section \ref{sec:Pf_MainThm}
is devoted to the proof of the main theorem.

\begin{table}[ht]
\small
\centering
\caption{
$95$ Reid orbifold $K3$ surfaces
in weighted project space
$\mathbb{CP}(\bs{w})$
($\bs{w}\in\mathbb{C}^{4}$)
and $\dim_{\mathbb{C}}\mathscr{M}^{*}$
}\label{table:main_result}
%\vspace{-10pt}
\begin{tabular}{cp{3mm}c}
\begin{tabular}{c>{\centering}p{10em}c}
\hline
\hline
Label
& $X_{d}\subset \mathbb{CP}(\bs{w})$
& $\dim_{\mathbb{C}}\mathscr{M}^*$ \\ 
\hline
\hline
1 & $X_{4}\subset\mathbb{CP}(1,1,1,1)$ & $90$ \\
2 & $X_{6}\subset\mathbb{CP}(1,1,1,3)$ & $90$ \\
%\hline %
3 & $X_{5}\subset\mathbb{CP}(1,1,1,2)$ & $84$\\
4 & $X_{12}\subset\mathbb{CP}(1,1,4,6)$ & $84$\\ 
%\hline %
5 & $X_{8}\subset\mathbb{CP}(1,1,2,4)$ & $78$\\ 
6 & $X_{10}\subset\mathbb{CP}(1,1,3,5)$ & $80$\\ 
%\hline %
7 & $X_{6}\subset\mathbb{CP}(1,1,2,2)$ & $72$\\ 
8 & $X_{7}\subset\mathbb{CP}(1,1,2,3)$ & $74$\\ 
9 & $X_{9}\subset\mathbb{CP}(1,1,3,4)$ & $76$\\ 
%\hline %
10 & $X_{10}\subset\mathbb{CP}(1,2,2,5)$ & $60$\\ 
11 & $X_{18}\subset\mathbb{CP}(1,2,6,9)$ & $62$\\ 
%\hline %
12 & $X_{8}\subset\mathbb{CP}(1,2,2,3)$ & $56$\\ 
13 & $X_{12}\subset\mathbb{CP}(1,2,3,6)$ & $58$\\ 
14 & $X_{14}\subset\mathbb{CP}(1,2,4,7)$ & $58$\\ 
15 & $X_{16}\subset\mathbb{CP}(1,2,5,8)$ & $60$\\ 
%\hline %
16 & $X_{9}\subset\mathbb{CP}(1,2,3,3)$ & $54$\\ 
17 & $X_{10}\subset\mathbb{CP}(1,2,3,4)$ & $54$\\ 
18 & $X_{11}\subset\mathbb{CP}(1,2,3,5)$ & $56$\\ 
19 & $X_{12}\subset\mathbb{CP}(1,2,4,5)$ & $54$\\ 
20 & $X_{15}\subset\mathbb{CP}(1,2,5,7)$ & $58$\\ 
21 & $X_{24}\subset\mathbb{CP}(1,3,8,12)$ & $56$\\ 
%\hline %
22 & $X_{16}\subset\mathbb{CP}(1,3,4,8)$ & $52$\\ 
23 & $X_{18}\subset\mathbb{CP}(1,3,5,9)$ & $52$\\ 
24 & $X_{22}\subset\mathbb{CP}(1,3,7,11)$ & $54$\\ 
25 & $X_{30}\subset\mathbb{CP}(1,4,10,15)$ & $52$\\ 
26 & $X_{12}\subset\mathbb{CP}(1,3,4,4)$ & $48$ \\ 
27 & $X_{13}\subset\mathbb{CP}(1,3,4,5)$ & $48$\\ 
28 & $X_{14}\subset\mathbb{CP}(2,2,3,7)$ & $38$\\ 
29 & $X_{15}\subset\mathbb{CP}(1,3,4,7)$ & $50$\\ 
30 & $X_{15}\subset\mathbb{CP}(1,3,5,6)$ & $48$ \\ 
31 & $X_{20}\subset\mathbb{CP}(1,4,5,10)$ & $48$ \\ 
32 & $X_{21}\subset\mathbb{CP}(1,3,7,10)$ & $52$\\ 
33 & $X_{22}\subset\mathbb{CP}(1,4,6,11)$ & $48$\\ 
34 & $X_{28}\subset\mathbb{CP}(1,4,9,14)$ & $50$\\ 
35 & $X_{36}\subset\mathbb{CP}(1,5,12,18)$ & $50$\\ 
36 & $X_{42}\subset\mathbb{CP}(1,6,14,21)$ & $48$\\ 
%\hline %
37 & $X_{12}\subset\mathbb{CP}(2,2,3,5)$ & $36$\\ 
38 & $X_{16}\subset\mathbb{CP}(1,4,5,6)$ & $44$ \\ 
39 & $X_{18}\subset\mathbb{CP}(1,4,6,7)$ & $44$ \\ 
40 & $X_{26}\subset\mathbb{CP}(1,5,7,13)$ & $46$\\ 
41 & $X_{30}\subset\mathbb{CP}(1,6,8,15)$ & $44$ \\ 
%\hline %
42 & $X_{12}\subset\mathbb{CP}(2,3,3,4)$ & $32$\\ 
43 & $X_{18}\subset\mathbb{CP}(2,3,4,9)$ & $32$\\ 
44 & $X_{21}\subset\mathbb{CP}(1,5,7,8)$ & $42$\\ 
45 & $X_{24}\subset\mathbb{CP}(1,6,8,9)$ & $40$\\ 
46 & $X_{30}\subset\mathbb{CP}(2,3,10,15)$ & $34$\\ 
%\hline %
47 & $X_{14}\subset\mathbb{CP}(2,3,4,5)$ & $30$\\ 
48 & $X_{16}\subset\mathbb{CP}(2,3,4,7)$ & $30$\\ 
%\hline %
%%
%$\vdots$ &&& \\
%??? & $13$ & $X_{42}\subset\mathbb{CP}^{2}(1,6,14,21)$ & $48$ \\
%$\vdots$ &&& \\
\hline
\hline
\end{tabular}
&&
\begin{tabular}{c>{\centering}p{10em}c}
\hline
\hline
Label
& $X_{d}\subset \mathbb{CP}(\bs{w})$
& $\dim_{\mathbb{C}}
\mathscr{M}^*$ \\ 
\hline
\hline
49 & $X_{20}\subset\mathbb{CP}(2,3,5,10)$ & $32$ \\ 
50 & $X_{22}\subset\mathbb{CP}(2,4,5,11)$ & $28$\\ 
51 & $X_{24}\subset\mathbb{CP}(2,3,7,12)$ & $32$ \\ 
52 & $X_{26}\subset\mathbb{CP}(2,3,8,13)$ & $32$\\ 
%\hline %
53 & $X_{15}\subset\mathbb{CP}(2,3,5,5)$ & $30$\\ 
54 & $X_{15}\subset\mathbb{CP}(3,3,4,5)$ & $26$\\ 
55 & $X_{17}\subset\mathbb{CP}(2,3,5,7)$ & $30$ \\ 
56 & $X_{18}\subset\mathbb{CP}(2,3,5,8)$ & $30$ \\ 
57 & $X_{20}\subset\mathbb{CP}(2,4,5,9)$ & $26$\\ 
58 & $X_{21}\subset\mathbb{CP}(2,3,7,9)$ & $30$ \\ 
59 & $X_{24}\subset\mathbb{CP}(2,3,8,11)$ & $30$ \\ 
60 & $X_{26}\subset\mathbb{CP}(2,5,6,13)$ & $26$\\ 
61 & $X_{30}\subset\mathbb{CP}(2,6,7,15)$ & $24$\\ 
%\hline %
62 & $X_{18}\subset\mathbb{CP}(3,4,5,6)$ & $22$\\ 
63 & $X_{20}\subset\mathbb{CP}(2,5,6,7)$ & $24$\\ 
64 & $X_{24}\subset\mathbb{CP}(3,4,5,12)$ & $32$ \\ 
65 & $X_{32}\subset\mathbb{CP}(2,5,9,16)$ & $26$ \\ 
66 & $X_{42}\subset\mathbb{CP}(3,4,14,21)$ & $24$\\ 
%\hline %
67 & $X_{19}\subset\mathbb{CP}(3,4,5,7)$ & $22$\\ 
68 & $X_{20}\subset\mathbb{CP}(3,4,5,8)$ & $22$\\ 
69 & $X_{21}\subset\mathbb{CP}(3,5,6,7)$ & $20$\\ 
70 & $X_{27}\subset\mathbb{CP}(2,5,9,11)$ & $24$ \\ 
71 & $X_{28}\subset\mathbb{CP}(3,4,7,14)$ & $22$ \\ 
72 & $X_{30}\subset\mathbb{CP}(4,5,6,15)$ & $18$\\ 
73 & $X_{34}\subset\mathbb{CP}(3,4,10,17)$ & $22$\\ 
74 & $X_{36}\subset\mathbb{CP}(3,4,11,18)$ & $22$\\ 
75 & $X_{48}\subset\mathbb{CP}(3,5,16,24)$ & $22$ \\ 
76 & $X_{24}\subset\mathbb{CP}(3,4,7,10)$ & $20$  \\ 
77 & $X_{24}\subset\mathbb{CP}(4,5,6,9)$ & $16$\\ 
78 & $X_{30}\subset\mathbb{CP}(3,4,10,13)$ & $20$\\ 
79 & $X_{32}\subset\mathbb{CP}(4,5,7,16)$ & $0$\\ 
80 & $X_{34}\subset\mathbb{CP}(4,6,7,17)$ & $16$\\ 
81 & $X_{38}\subset\mathbb{CP}(3,5,11,19)$ & $20$\\ 
82 & $X_{54}\subset\mathbb{CP}(4,5,18,27)$ & $18$ \\ 
%\hline %
83 & $X_{25}\subset\mathbb{CP}(4,5,7,9)$ & $16$\\ 
84 & $X_{27}\subset\mathbb{CP}(5,6,7,9)$ & $14$\\ 
85 & $X_{28}\subset\mathbb{CP}(4,6,7,11)$ & $14$\\ 
86 & $X_{33}\subset\mathbb{CP}(3,5,11,14)$ & $18$\\ 
87 & $X_{38}\subset\mathbb{CP}(5,6,8,19)$ & $14$\\ 
88 & $X_{40}\subset\mathbb{CP}(5,7,8,20)$ & $14$\\ 
89 & $X_{42}\subset\mathbb{CP}(2,5,14,21)$ & $28$\\ 
90 & $X_{44}\subset\mathbb{CP}(4,5,13,22)$ & $16$ \\ 
91 & $X_{66}\subset\mathbb{CP}(5,6,22,33)$ & $14$ \\ 
%\hline%
92 & $X_{24}\subset\mathbb{CP}(3,6,7,8)$ & $18$ \\ 
93 & $X_{30}\subset\mathbb{CP}(5,6,8,11)$ & $12$\\ 
94 & $X_{36}\subset\mathbb{CP}(7,8,9,12)$ & $10$ \\ 
95 & $X_{50}\subset\mathbb{CP}(7,8,10,25)$ & $10$ \\ 
\hline
\hline
\rule{0pt}{8pt}\\
\end{tabular}
\end{tabular}
\end{table}

\clearpage

\section{Preliminaries}\label{sec:Pre}

In this section,
we recall the basic definitions and conventions used
throughout the paper. 
We adopt the standard terminology of Sasakian geometry
and contact instantons as in \cite{BG, BH}.

\subsection{Sasakian manifold}

We begin by recalling the notion of an almost contact metric structure on an odd-dimensional Riemannian manifold.
Let $M$ be a $(2n+1)$-dimensional smooth manifold
and $g$ be a Riemannian metric on $M$
with Levi-Civita connection $\nabla^{g}$.

\begin{dfn}
An \textit{almost contact metric structure}
$(\xi, \eta, \Phi, g)$
on $M$ consists of a vector field $\xi$,
a $1$-form $\eta$,
an endomorphism $\Phi:TM \rightarrow TM$
and the metric $g$ satisfying the following conditions:
\begin{align}
  &\eta(\xi)=1,\qquad \Phi^2=-\mathrm{Id}+\eta\otimes\xi,\label{eqn:contact_str}\\
  &g(\Phi X,\Phi Y)=g(X,Y)-\eta(X)\eta(Y),\label{eqn:contact_str2}
\end{align}
for all vector fields $X, Y$ on $M$.
\end{dfn}

Let $(\xi, \eta, \Phi, g)$ be an almost contact metric structure on $M$.
The first equation of
\eqref{eqn:contact_str}
implies that $\xi$ is nowhere vanishing on $M$.
Then $\xi$ generates a $1$-dimensional subbundle of $TM$,
which we write $L_{\xi}$.
Denote by $\mathcal{F}_{\xi}$
the $1$-dimensional foliation
generated by $\xi$,
which is called the \textit{characteristic foliation} of $M$.
The kernel $H:=\mathrm{Ker}(\eta)(\subset TM)$
of $\eta$ is called the \textit{contact distribution}
for $M$.
Then $TM$ is decomposed into $TM=H\oplus L_{\xi}$.
It follows from \cite[Theorem 4.1]{Blair}
that $\Phi\xi=0$ and $\eta\circ\Phi=0$ hold
and that the rank of $\Phi$ is equal to $2n$.
By substituting $Y=\xi$ into \eqref{eqn:contact_str2},
we get $g(X,\xi)=\eta(X)$.
In particular, $\xi$ is a unit vector field on $M$
and $g(H,L_{\xi})=\{0\}$.
The normal bundle $TM/L_{\xi}$ is identified with $H$.
The metric $g^{T}:=g|_{H}$
on $H$ is called the \textit{transverse metric}.
We write $H^{*}$ as the dual of $H$
and $\pi_{H}:TM\to H$ as the orthogonal projection with respect to $TM=H\oplus L_{\xi}$.

The restriction $J^{T}:=\Phi|_{H}$ of $\Phi$ to
$H$ gives an almost complex structure on $H$,
which is called a transverse almost complex structure.
By \eqref{eqn:contact_str2},
the transverse metric $g^{T}$ 
on $H$ is Hermitian with respect to $J^{T}$.
We define the fundamental $2$-form $\omega$
on $M$ by
\begin{equation}\label{eqn:dfn_omega}
\omega(X,Y):=g(X,\Phi Y).
\end{equation}

\begin{dfn}
A almost contact metric structure
$(\xi, \eta, \Phi, g)$ on $M$ is called a \textit{contact metric structure}
if $d\eta=\omega$ holds.
Then $M$ equipped with this structure is called
a \textit{contact metric manifold}.
\end{dfn}

For any contact metric structure
$(\xi, \eta, \Phi, g)$
on $M$,
we have
$d\eta(\xi,\cdot)=0$
and $\eta\wedge (d\eta)^{n}\neq 0$,
from which $\eta$ is a contact form
on $M$ and $\xi$ is the Reeb vector field
associated to $\eta$.

\begin{dfn}
A contact metric structure
$(\xi, \eta, \Phi,g)$ on $M$
is called a \textit{Sasakian structure}
if it satisfies the following condition:
\begin{equation}\label{eqn:eqn:sasaki}
(\nabla^{g}_{X}\Phi)Y=g(X,Y)\xi-\eta(Y)X,
\end{equation}
for all vector fields $X,Y$ on $M$.
Then $M$ equipped with this structure
is called a
\textit{Sasakian manifold}.
\end{dfn}

Let $(M, \xi, \eta, \Phi,g)$
be a Sasakian manifold.
The condition \eqref{eqn:eqn:sasaki}
implies that $\nabla^{g}_{X}\xi=-\Phi X$,
from which we get:
\begin{equation}
(\mathcal{L}_{\xi}g)(X,Y)=g(\nabla^{g}_{X}\xi,Y)+g(X,\nabla^{g}_{Y}\xi)=0,
\end{equation}
where $\mathcal{L}_{\xi}$ denotes the Lie derivative of $\xi$.
Hence $\xi$ is a Killing vector field on $M$.
It follows from
\cite[Proposition 6.4.8]{BG} that
the characteristic foliation $\mathcal{F}_{\xi}$
is a Riemannian foliation.

The condition \eqref{eqn:eqn:sasaki}
implies that the Nijenhuis tensor of $J^{T}$
on $H$ vanishes,
that is, $J^{T}$ is integrable.
The restriction $\omega^{T}:=\omega|_{H}$
of the fundamental $2$-form $\omega$ to $H$
is called the \textit{transverse K\"ahler $2$-form}
associated with $(J^{T},g^{T})$.
The triplet $(J^{T},g^{T},\omega^{T})$
is called a transverse K\"ahler structure on $H$.

For any linear connection $\nabla$ on $M$,
the \textit{transverse connection} $\nabla^{T}$
on $H$ associated with $\nabla$
is defined by
\begin{equation}
\nabla^{T}_{X}Y = \pi_{H}(\nabla_{X}Y),\quad
\nabla^{T}_{V}Y=\pi_{H}([V,Y]),
\end{equation}
for all $X,Y\in\Gamma(M,H)$ and $V\in\Gamma(M,L_{\xi})$
(\cite[(2.5.7)]{BG}).
Here,
the torsion
$T^{T}\in\Gamma(M,\Lambda^{2}H^{*}\otimes H)$
of $\nabla^{T}$
is defined by
$T^{T}(X,Y)=
\nabla^{T}_{X}Y-\nabla^{T}_{Y}X-\pi_{H}([X,Y])$
for $X,Y\in\Gamma(M,H)$.
The transverse connection
associated with $\nabla^{g}$
is called the \textit{transverse Levi-Civita connection}
for $M$,
which we write $\nabla^{*}$ instead of $(\nabla^{g})^{T}$.
It is known that
$\nabla^{*}$
is a unique transverse
connection on $H$ which is torsion-free and compatible with $g^{T}$.

\subsection{Transverse Calabi-Yau structures}

We begin by recalling the notion of
basic cohomology.
Let $M$ be a compact $5$-dimensional Sasakian manifold 
with Sasakian structure $(\xi, \eta, \Phi,g)$.
A $k$-form $\alpha$
is said to be \textit{basic}
if it satisfies 
$\iota_{\xi}\alpha=0$
and
$\mathcal{L}_{\xi}\alpha=0$,
where $\iota_{\xi}$
denotes the interior product of $\xi$.
Let $\Omega^{k}_{B}(M)$
denote the space of basic $k$-forms on $M$.
It is verified that
$d\alpha$
is basic for any $\alpha\in\Omega^{k}_{B}(M)$.
Hence the exterior derivative $d$ maps
$\Omega^{k}_{B}(M)$ into $\Omega^{k+1}_{B}(M)$.
Then the \textit{basic cohomology}
$H^{k}_{B}(M)$ is defined by
\begin{equation}
H^{k}_{B}(M)
=\mathrm{Ker}(d:\Omega^{k}_{B}(M)\to\Omega^{k+1}_{B}(M))/
\mathrm{Im}(d:\Omega^{k-1}_{B}(M)\to\Omega^{k}_{B}(M)).
\end{equation}

The contact distribution $H$ has
the structure of
a complex vector bundle over $M$
associated with $J^{T}$.
The complexification $H\otimes\mathbb{C}$
of $H$ admits a decomposition into eigenspaces of $J^{T}$:
\begin{equation}
H\otimes\mathbb{C}=H^{(1,0)}\oplus H^{(0,1)},
\end{equation}
where $H^{(1,0)}:=\{X\in H\otimes\mathbb{C}\mid J^{T}X=iX\}$
and $H^{(0,1)}:=\{X\in H\otimes\mathbb{C}\mid J^{T}X=-iX\}$.
Let
$\Omega^{k}_{B}(M)\otimes\mathbb{C}$
denote the space of $\mathbb{C}$-valued basic $k$-forms
and $\Omega^{(p,q)}_{B}(H)$
denote the space of forms of type $(p,q)$,
that is,
\begin{equation}
\Omega^{(p,q)}_{B}(M)
:=\{\alpha\in\Omega^{p+q}_{B}(M)\otimes\mathbb{C}
\mid
\alpha \in
\Lambda^{p}(H^{(1,0)})^{*}
\otimes
\Lambda^{q}(H^{(0,1)})^{*}
\}.
\end{equation}
The elements of
$\Omega^{(p,q)}_{B}(M)$
are called \textit{basic $(p,q)$-forms}.
Then
$\Omega^{k}_{B}(M)\otimes\mathbb{C}$
is decomposed as follows:
\begin{equation}
\Omega^{k}_{B}(M)\otimes\mathbb{C}
=\bigoplus_{p+q=k}\Omega^{p,q}_{B}(M).
\end{equation}
Let $\nabla^{T}$ be
a transverse connection on $H$
and $\Omega^{T}$ denote the transverse 
curvature $2$-form of $\nabla^{T}$.
It follows from \cite[Proposition 7.5.16]{BG}
that $\Omega^{T}$ is a basic $(1,1)$-form.
As in ordinary Chern-Weil theory,
the first Chern form $c_{1}(H,\nabla^{T})$
is defined by
\begin{equation}
c_{1}(H,\nabla^{T})=-\dfrac{1}{2\pi i}\operatorname{tr}\Omega^{T}.
\end{equation}
Then $c_{1}(H,\nabla^{T})$
is a basic $(1,1)$-form.
Furthermore, the basic cohomology class
of $c_{1}(H,\nabla^{T})$ is independent of the choice of
transverse connections on $H$,
which is called the basic first Chern class
of $M$ (\cite[Theorem/Definition 7.5.17]{BG}).
We denote by $c_{1}^{B}(H)$ this basic cohomology class
and call it
the \textit{basic first Chern class}
of $M$.
It should be noted that
$c_{1}^{B}(H)$
is an invariant of the underlying contact structure.

\begin{dfn}
A compact $5$-dimensional Sasakian manifold
satisfying $c_{1}^{B}(H)=0$
is called a Sasakian manifold with a \textit{transverse Calabi-Yau structure}.
\end{dfn}

Let $M$ be a compact $5$-dimensional Sasakian manifold
with a transverse Calabi-Yau structure.
It follows from \cite[Theorem 8.1.14]{BG}
that $M$ is quasi-regular in the sense of \cite[Definition 6.1.25]{BG}.
Hence, by \cite[Theorem 7.1.13]{BG},
the leaf space $X:=M/\mathcal{F}_{\xi}$
is a compact K\"ahler orbifold
and the canonical projection
$\pi: M\to X$
is an orbifold Riemannian submersion.
Moreover, since 
$c_{1}^{B}(M)=0$,
it follows from
\cite[Proposition 7.5.23]{BG}
that 
the orbifold first Chern class of $M/\mathcal{F}_{\xi}$
vanishes.
A compact K\"ahler orbifold with
vanishing the orbifold first Chern class
is called a \textit{Calabi-Yau orbifold}.

\subsection{Contact instantons and their moduli spaces}
Let $M$ be a $5$-dimensional Sasakian manifold
with Sasakian structure $(\xi,\eta,\Phi,g)$.
A $k$-form $\alpha$ on $M$ is said to be \textit{transverse}
if $\iota_{\xi}\alpha = 0$ holds.
Let
$\Omega^{k}_{H}(M) := \Gamma(M, \Lambda^{k}H^{*})$
denote the space of transverse $k$-forms on $M$.
Then we have
\begin{equation}\label{eqn:O2_O2H_etaO1H}
\Omega^{2}(M)=\Omega^{2}_{H}(M)
\oplus (\eta\wedge\Omega^{1}_{H}(M)).
\end{equation}
The \textit{transverse Hodge star operator}
$\ast_{T} \colon \Omega^{k}_{H}(M)
  \to \Omega^{4-k}_{H}(M)$
is defined as follows
(\cite[(7.2.2)]{BG}):
\begin{equation}\label{eqn:transverse-Hodge}
  \ast_{T}\alpha := \ast(\eta \wedge \alpha),
\end{equation}
where $\ast$ denotes the Hodge star operator on $M$.
Since
$\ast_{T}^{2} = \mathrm{Id}$
holds on $\Omega^{2}_{H}(M)$,
we have the eigenspace
decomposition of $\Omega^{2}_{H}(M)$
for $\ast_{T}$:
\begin{equation}\label{eqn:SD-ASD-decomp}
  \Omega^{2}_{H}(M)
  = \Omega^{+}_{H}(M) \oplus \Omega^{-}_{H}(M).
\end{equation}
A transverse $2$-form $\alpha$ is said to be
\textit{self-dual} (resp.\ \textit{anti-self-dual})
if $\ast_{T}\alpha = \alpha$
(resp.\ $\ast_{T}\alpha = -\alpha$).
The fundamental $2$-form $\omega$
is a transverse $2$-form and is self-dual
in this sense.

Now let $G$ be a compact connected Lie group,
$P$ be a principal $G$-bundle over $M$,
and $\mathfrak{g}_{P} := P \times_{\mathrm{Ad}} \mathfrak{g}$
denote the adjoint bundle.
Denote by
$F^{\nabla} \in \Omega^{2}(M)\otimes
\mathfrak{g}_{P}$
the curvature $2$-form of
a connection
$\nabla$ on
$P$.
Since
$\Omega^{2}(M)\otimes
\mathfrak{g}_{P}=
(\Omega^{2}_{H}(M)\otimes\mathfrak{g}_{P})
\oplus
((\eta\wedge
\Omega^{1}_{H}(M))\otimes\mathfrak{g}_{P})
$,
we write 
$F^{\nabla}$
as $F^{\nabla} = F^{\nabla}_{H} + \eta \wedge \alpha$,
where $F^{\nabla}_{H} \in \Omega^{2}_{H}(M, \mathfrak{g}_{P})$
is the transverse component of $F^{\nabla}$
and $\alpha$ is in $\Omega^{1}_{H}(M, \mathfrak{g}_{P})$.

\begin{dfn}\label{def:contact-instanton}
A connection
$\nabla$ on $P$ is called
a \textit{self-dual contact instanton}
(SD contact instanton)
if $F^{\nabla}_{H}$ is self-dual, i.e.,
\begin{equation}\label{def:SD}
  \ast_{T} F^{\nabla}_{H} = F^{\nabla}_{H},
\end{equation}
and an \textit{anti-self-dual contact instanton}
(ASD contact instanton)
if $F^{\nabla}_{H}$ is anti-self-dual, i.e.,
\begin{equation}\label{def:ASD}
\ast_{T} F^{\nabla}_{H} = -F^{\nabla}_{H}.
\end{equation}
\end{dfn}

Let $\nabla$ be a connection on $P$.
The conditions 
\eqref{def:SD}
and \eqref{def:ASD}
are equivalent
to
\begin{equation}
\ast F^{\nabla} = \eta \wedge F^{\nabla},
\quad
\ast F^{\nabla} = -\eta \wedge F^{\nabla},
\end{equation}
respectively.
It is known that
any ASD contact instanton $\nabla$
on $P$
satisfies
$d^{\nabla}\ast F^{\nabla}=0$,
that is, $\nabla$
is a Yang-Mills connection on $P$.
On the other hand,
SD contact instantons
are not necessarily Yang-Mills connections.
The group $\mathscr{G}=\Gamma(M,P\times_{\operatorname{Ad}}G)$ of
gauge transformations of $P$
acts on the space of connections
on $P$.
If $\nabla$
is a ASD contact instanton,
then $u^{*}\nabla$
is also an ASD contact instanton
for all $u\in\mathscr{G}$.

The center $Z(G)$ of $G$
is naturally embedded in $\mathscr{G}$
as constant maps.
For any connection $\nabla$
on $P$,
the isotropy subgroup
$\mathscr{G}_{\nabla}$
of $\mathscr{G}$
at $\nabla$
contains $Z(G)$.
A connection $\nabla$
on $P$ is said to be \textit{irreducible}
if $\mathscr{G}_{\nabla}=Z(G)$ holds.

\begin{dfn}\label{def:moduli}
Let $\mathscr{C}^{*}_{\mathrm{ASD}}$
denote the space of irreducible ASD contact instantons on $P$.
The quotient
\begin{equation}
\mathscr{M}^{*} := \mathscr{C}^{*}_{\mathrm{ASD}} / \mathscr{G}
\end{equation}
is called the
\textit{moduli space of irreducible ASD contact instantons}
on $P$.
\end{dfn}
In the case when
the Sasakian manifold $M$
has a transverse Calabi-Yau structure,
it follows from
\cite[Corollary 4.4, Proposition 4.19]{BH}
that the moduli space $\mathscr{M}^{*}$
is a smooth hyperk\"ahler manifold,
possibly empty.

\section{Proof of Main theorem}\label{sec:Pf_MainThm}

Let $(M,\xi,\eta,\Phi,g)$
be a compact simply-connected,
$5$-dimensional Sasakian manifold
with a transverse Calabi-Yau structure.
The transverse Levi-Civita connection
$\nabla^{*} := (\nabla^{g})^{T}$
is a connection on the $\mathrm{SO}(4)$-frame bundle of $H$.
The decomposition \eqref{eqn:SD-ASD-decomp}
at the bundle level gives
\begin{equation}\label{eqn:bundle-decomp}
  \Lambda^{2}H^{*}
  = \Lambda^{+}H^{*} \oplus \Lambda^{-}H^{*}.
\end{equation}
The Lie algebra of the structure group $\mathrm{SO}(4)$
is decomposed into
$\mathfrak{so}(4) \cong \mathfrak{so}(3)_{+}
  \oplus \mathfrak{so}(3)_{-}$.
The anti-self-dual component of $\nabla^{*}$
defines a connection on the $\mathrm{SO}(3)$-bundle
whose adjoint bundle is
$\mathfrak{g}_{P} = \Lambda^{-}H^{*}$.
Since $M$ has a transverse Calabi-Yau structure,
the transverse Ricci curvature vanishes,
and hence by \cite[Proposition 2.11]{BH}
this connection is an ASD contact instanton.

Reid (\cite{Reid1, Reid2})
classified
orbifold $K3$ surfaces
among quasi-smooth, well-formed hypersurface
in the weighted projective space
$\mathbb{CP}(w_{0},w_{1},w_{2},w_{3})$.
His classification list,
which consists of $95$ families of orbifold $K3$ surfaces,
provides a systematic source of the leaf spaces
of Sasakian manifolds with transverse Calabi-Yau structures.
We call such an orbifold $K3$ surface
a \textit{Reid orbifold $K3$ surface}.
Each family in Reid's list is specified by
a weight vector
$(w_{0},w_{1},w_{2},w_{3})$
with
$w_{0} \leq w_{1} \leq w_{2} \leq w_{3}$
and
$\gcd(w_{i},w_{j},w_{k}) = 1$
for every triple of distinct indices
$\{i,j,k\} \subset \{0,1,2,3\}$.
The corresponding hypersurface
$X_{d} \subset
\mathbb{CP}(w_{0},w_{1},w_{2},w_{3})$
has degree $d = w_{0}+w_{1}+w_{2}+w_{3}$,
so that $X_{d}$ is a Calabi-Yau orbifold.
Since $X_{d}$ is quasi-smooth,
its only singularities are
cyclic quotient singularities
inherited from $\mathbb{CP}(w_{0},w_{1},w_{2},w_{3})$.
These are du~Val singularities
of type $A_{m_{j}-1}$,
where $m_{j}$ is the order
of the isotropy group
at the $j$-th singular point
for $j=1,\dotsc,k$.
Fletcher (\cite{Fletcher}) determined
the singularity data
$\{m_{1},\dotsc,m_{k}\}$
for all $95$ families.

In what follows, we first prove the irreducibility of the transverse
Levi-Civita connection $\nabla^{*}$ on $H$ and the moduli space
$\mathscr{M}^{*}$
by using the singularity
data of $X=M/\mathcal{F}_{\xi}$.

\begin{thm}[{\cite[(5.9)]{BH}}]\label{thm:irreducibility}
Suppose that the leaf space $X$
is a Reid orbifold $K3$ surface
with singularity data $\{m_{1},\dotsc,m_{k}\}$.
Then,
$\nabla^{*}$ is irreducible
if the following relation holds:
\begin{equation}\label{eq:irred-condition}
\sum_{j=1}^{k}\frac{m_{j}^{2}-1}{m_{j}} \neq 24,
\end{equation}
\end{thm}

Motivated by the above theorem,
we call the left hand side of
\eqref{eq:irred-condition}
the \textit{degree}
of $\nabla^{*}$, 
which we write $\mathrm{deg}(\nabla^{*})$:
\begin{equation}\label{eqn:dfn_deg}
\mathrm{deg}(\nabla^{*}):=\sum_{j=1}^{k}\frac{m_{j}^{2}-1}{m_{j}}.
\end{equation}
In order to show
that $\nabla^{*}$ is irreducible
we calculate the value of
$\mathrm{deg}(\nabla^{*})$ for each
Sasakian manifolds with transverse Calabi-Yau structures
whose leaf spaces
are Reid orbifold $K3$ surface $X_{d}\subset\mathbb{CP}(w_{0},w_{1},w_{2},w_{3})$
as listed in Table \ref{table:main_result}.

\begin{ex}[{Label 33 in Table \ref{table:main_result}}]\label{ex:label33}
Let $X$ be a Reid orbifold $K3$ surface of degree $22$ in the weighted projective space $\mathbb{CP}(1,4,6,11)$.
It follows from \cite{Fletcher} that
$X$ has three singularities
whose types are $A_{3},A_{1},A_{5}$. 
By substituting $m_{1}=4, m_{2}=2$
and $m_{3}=6$ into the right hand side of 
\eqref{eqn:dfn_deg}, we get
\begin{equation}
\mathrm{deg}(\nabla^{*})
=\dfrac{4^{2}-1}{4}+\dfrac{2^{2}-1}{2}+\dfrac{6^{2}-1}{6}
=\frac{133}{12}\neq 24.
\end{equation}
Thus,
by Theorem \ref{thm:irreducibility},
the transverse Levi-Civita connection $\nabla^{*}$
is irreducible.
\end{ex}

Table \ref{table:degree} below shows the result
of determining $\mathrm{deg}(\nabla^{*})$ 
for each
Sasakian manifolds with transverse Calabi-Yau structures
whose leaf spaces
are Reid orbifold $K3$ surface $X_{d}\subset\mathbb{CP}(\bs{w})$
as listed in Table \ref{table:main_result}.
The labeling convention in Table \ref{table:degree} is the same as that used for the 
Reid orbifold $K3$ surfaces
in $\mathbb{CP}(\bs{w})$ in Table \ref{table:main_result}.
Then we find that $\mathrm{deg}(\nabla^{*})$
is not equal to $24$,
so that $\nabla^{*}$ is irreducible
due to Theorem \ref{thm:irreducibility}.
Hence the moduli space $\mathscr{M}^{*}$
becomes a non-empty hyperk\"ahler manifold.

Next, we determine the complex dimension of $\mathscr{M}^{*}$
by using the following theorem.

\begin{thm}[{\cite[p.\,591]{BH}}]\label{thm:dim}
Under the same setting of Theorem \ref{thm:irreducibility}, we have
\begin{equation}
\dim_{\mathbb{C}}\mathscr{M}^{*}=90-2\sum_{j=1}^{k}(2m_{j}-1).
\end{equation}
\end{thm}

\begin{ex}[continued in Example \ref{ex:label33}]
By using Theorem \ref{thm:dim},
the dimension of the moduli space is computed as
$\dim_{\mathbb{C}}\mathscr{M}^{*}=90-2(7+3+11)=48$.
We conclude that
$\mathscr{M}^{*}$ is a hyperk\"ahler manifold of complex dimension
$48$.
\end{ex}

A similar calculation
shows
the result
of determining $\dim_{\mathbb{C}}\mathscr{M}^{*}$ 
as in Table \ref{table:main_result}.
Thus, we have completed the proof of Main theorem
stated in Introduction.

\begin{table}[ht]
\centering
\small
\caption{
The singularity data of 
Reid orbifold $K3$ surfaces $X_{d}$
in $\mathbb{CP}(\bs{w})$
and $\mathrm{deg}(\nabla^{*})$}\label{table:degree}
\begin{tabular}{cp{3mm}c}
\begin{tabular}{ccc}
\hline
\hline
Label & Singularities of $X_{d}$ & $\mathrm{deg}(\nabla^{*})$ \\ 
\hline
\hline
1 & no singularities & $0$\\
2 & no singularities & $0$\\
%\hline %
3 & $A_{1}$ & $3/2$\\
4 & $A_{1}$ &  $3/2$\\ 
%\hline %
5 & $2\times A_{1}$ & $3$\\ 
6 & $A_{2}$ & $8/3$\\ 
%\hline %
7 & $3\times A_{1}$ & $9/2$\\ 
8 & $A_{1}$, $A_{2}$ & $25/6$\\ 
9 & $A_{3}$ & $15/4$\\ 
%\hline %
10& $5\times A_{1}$ & $15/2$\\ 
11& $3\times A_{1}$, $A_{2}$ & $43/6$\\ 
%\hline %
12& $4\times A_{1}$, $A_{2}$ & $26/3$\\ 
13& $2\times A_{1}$, $2\times A_{2}$ & $25/3$\\ 
14& $3\times A_{1}$, $A_{3}$ & $33/4$\\ 
15& $2\times A_{1}$, $A_{4}$ & $39/5$\\ 
%\hline %
16& $A_{1}$, $3\times A_{2}$ & $19/2$\\ 
17& $2\times A_{1}$, $A_{2}$, $A_{3}$ & $113/12$\\ 
18& $A_{1}$, $A_{2}$, $A_{4}$ & $269/30$\\ 
19& $3\times A_{1}$, $A_{4}$ & $93/10$\\ 
20& $A_{1}$, $A_{6}$ & $117/14$\\ 
\hline
\hline
\end{tabular}
&&
\begin{tabular}{ccc}
\hline
\hline
Label & Singularities of $X_{d}$ & $\mathrm{deg}(\nabla^{*})$ \\ 
\hline
\hline
21& $2\times A_{2}$, $A_{3}$ & $109/12$\\ 
%\hline %
22& $A_{2}$, $2\times A_{3}$ & $61/6$\\ 
23& $2\times A_{2}$, $A_{4}$ & $152/15$\\ 
24& $A_{2}$, $A_{6}$ & $200/21$\\ 
25& $A_{3}$, $A_{1}$, $A_{4}$ & $201/20$\\ 
26& $3\times A_{3}$ & $45/4$\\ 
27& $A_{2}$, $A_{3}$, $A_{4}$ &$673/60$\\ 
28& $7\times A_{1}$, $A_{2}$ &$79/6$\\ 
29& $A_{3}$, $A_{6}$ &$297/28$\\ 
30& $2\times A_{2}$, $A_{5}$ &$67/6$\\ 
31& $A_{1}$, $2\times A_{4}$ &$110/10$\\ 
32& $A_{9}$ &$99/10$\\ 
33 & $A_{3}$, $A_{1}$, $A_{5}$ &$133/12$ \\ 
34 & $A_{1}$, $A_{8}$ &$187/18$ \\ 
35 & $A_{4}$, $A_{5}$ &$319/30$ \\ 
36 & $A_{1}$, $A_{2}$, $A_{6}$ &$463/42$ \\ 
%\hline %
37 & $6\times A_{1}$, $A_{4}$ &$69/5$ \\ 
38 & $A_{1}$, $A_{4}$, $A_{5}$ &$182/15$ \\ 
39 & $A_{3}$, $A_{1}$, $A_{6}$ &$339/28$ \\ 
40 & $A_{4}$, $A_{6}$ &$408/35$ \\ 
\hline
\hline
\end{tabular}
\end{tabular}
\end{table}

\addtocounter{table}{-1}
\begin{table}[ht]
\centering
\small
\caption{(continued)}
\begin{tabular}{cp{3mm}c}
\begin{tabular}{ccc}
\hline
\hline
Label & Singularities of $X_{d}$ & $\mathrm{deg}(\nabla^{*})$ \\ 
\hline
\hline
41 & $A_{1}$, $A_{2}$, $A_{7}$ &$289/24$ \\ 
%\hline %
42 & $3\times A_{1}$, $4\times A_{2}$ &$91/6$ \\ 
43 & $4\times A_{1}$, $2\times A_{2}$, $A_{3}$ &$181/12$ \\ 
44 & $A_{4}$, $A_{7}$ &$507/40$ \\ 
45 & $A_{1}$, $A_{2}$, $A_{8}$ &$235/18$ \\ 
46 & $3\times A_{1}$, $2\times A_{2}$, $A_{4}$ &$439/30$ \\ 
%\hline %
47 & $3\times A_{1}$, $A_{2}$, $A_{3}$, $A_{4}$  &$943/60$ \\ 
48 & $4\times A_{1}$, $A_{2}$, $A_{6}$ &$326/21$ \\ 
49 & $2\times A_{1}$, $A_{2}$, $2\times A_{4}$ &$229/15$ \\ 
50 & $5\times A_{1}$, $A_{3}$, $A_{4}$ &$321/20$ \\ 
51 & $2\times A_{1}$, $2\times A_{2}$, $A_{6}$ &$319/21$ \\ 
52 & $3\times A_{1}$, $A_{2}$, $A_{7}$ &$361/24$ \\ 
%\hline %
53 & $A_{1}$, $3\times A_{4}$ &$159/10$ \\ 
54 & $5\times A_{2}$, $A_{3}$ &$205/12$ \\ 
55 & $A_{1}$, $A_{2}$, $A_{4}$, $A_{6}$ &$3323/210$ \\ 
56 & $2\times A_{1}$, $A_{4}$, $A_{7}$ &$627/40$ \\ 
57 & $5\times A_{1}$, $A_{8}$ &$295/18$ \\ 
58 & $A_{1}$, $2\times A_{2}$, $A_{8}$ &$283/18$ \\ 
59 & $3\times A_{1}$, $A_{10}$ &$339/22$ \\ 
60 & $4\times A_{1}$, $A_{4}$, $A_{5}$ &$499/30$ \\ 
61 & $5\times A_{1}$, $A_{2}$, $A_{6}$ &$715/42$ \\ 
62 & $3\times A_{2}$, $A_{3}$, $A_{1}$, $A_{4}$  &$361/20$ \\ 
63 & $3\times A_{1}$, $A_{5}$, $A_{6}$ &$361/21$ \\ 
64 & $2\times A_{1}$, $2\times A_{3}$, $A_{4}$ &$153/10$ \\ 
65 & $2\times A_{1}$, $A_{4}$, $A_{8}$ &$751/45$ \\ 
66 & $2\times A_{2}$, $A_{3}$, $A_{1}$, $A_{6}$ &$1465/84$ \\ 
%\hline %
67 & $A_{2}$, $A_{3}$, $A_{4}$, $A_{6}$ &$7591/420$ \\ 
68 & $A_{2}$, $2\times A_{3}$, $A_{7}$ &$433/24$ \\ 
%\hline %
\hline
\hline
\end{tabular}
&&
\begin{tabular}{c>{\centering}p{10em}c}
\hline
\hline
Label & Singularities of $X_{d}$ & $\mathrm{deg}(\nabla^{*})$ \\ 
\hline
\hline
69 & $3\times A_{2}$, $A_{4}$, $A_{5}$ &$559/30$ \\ 
70 & $A_{1}$, $A_{4}$, $A_{10}$ &$1893/110$ \\ 
71 & $A_{2}$, $A_{1}$, $2\times A_{6}$ &$751/42$ \\ 
72 & $A_{3}$, $2\times A_{1}$, $2\times A_{4}$, $A_{2}$ &$1141/60$ \\ 
73 & $A_{2}$, $A_{3}$, $A_{1}$, $A_{9}$ &$1069/60$ \\ 
74 & $2\times A_{2}$, $A_{1}$, $A_{10}$ &$1171/66$ \\ 
75 & $2\times A_{2}$, $A_{4}$, $A_{7}$&$2161/120$ \\ 
76 & $A_{1}$, $A_{6}$, $A_{9}$ &$639/35$ \\ 
77 & $2\times A_{1}$, $A_{4}$, $A_{2}$, $A_{8}$ &$871/45$ \\ 
78 & $A_{3}$, $A_{1}$, $A_{12}$ &$945/52$ \\ 
79 & $2\times A_{3}$, $2\times A_{4}$, $A_{6}$ &$1677/70$ \\ 
80 & $A_{3}$, $2\times A_{1}$, $A_{5}$, $A_{6}$ &$1633/84$ \\ 
81 & $A_{2}$, $A_{4}$, $A_{10}$ &$3032/165$ \\ 
82 & $A_{3}$, $A_{1}$, $A_{4}$, $A_{8}$ &$3409/180$ \\ 
%\hline %
83 & $A_{3}$, $A_{6}$, $A_{8}$ &$4913/252$ \\ 
84 & $A_{4}$, $A_{5}$, $A_{2}$, $A_{6}$ &$1411/210$ \\ 
85 & $2\times A_{1}$, $A_{5}$, $A_{10}$ &$1303/66$ \\ 
86 & $A_{4}$, $A_{13}$ &$1311/70$ \\ 
87 & $A_{4}$, $A_{5}$, $A_{1}$, $A_{7}$ &$2401/120$ \\ 
88 & $2\times A_{4}$, $A_{6}$, $A_{3}$ &$2829/140$ \\ 
89 & $3\times A_{1}$, $A_{4}$, $A_{6}$ &$1131/70$ \\ 
90 & $A_{1}$, $A_{4}$, $A_{12}$ &$2499/130$ \\ 
91 & $A_{4}$, $A_{1}$, $A_{2}$, $A_{10}$ &$6559/330$ \\ 
%\hline%
92 & $4\times A_{2}$, $A_{1}$, $A_{6}$ &$799/42$ \\ 
93 & $A_{1}$, $A_{7}$, $A_{10}$ &$1785/88$ \\ 
94 & $A_{6}$, $A_{7}$, $A_{3}$, $A_{2}$ &$3553/168$ \\ 
95 & $A_{6}$, $A_{7}$, $A_{1}$, $A_{4}$ &$5889/280$ \\ 
\hline
\hline
\rule{0pt}{7.5pt}\\
\end{tabular}
\end{tabular}
\end{table}

%\clearpage

\end{document}